\documentclass[12pt]{amsart}
\usepackage[margin=1.3in]{geometry}
\usepackage[all,cmtip]{xy}

\newtheorem*{theorem}{Theorem}
\newtheorem*{lemma}{Lemma}
\newtheorem*{claim}{Claim}
\newtheorem*{maintheorem}{Main Theorem}

\newcommand{\spec}{\operatorname{Spec}}
\newcommand{\height}{\operatorname{ht}}

\begin{document}

\title{Multiplicity of jet schemes of monomial schemes}

\author[C. Yuen]{Cornelia Yuen}
\address{Department of Mathematics, University of Michigan,
Ann Arbor, MI 48109, USA}
\email{{\tt oyuen@umich.edu}}

\begin{abstract}
    This article studies jet schemes of monomial schemes.  They are
    known to be equidimensional but usually are not reduced
    \cite{GowardSmith}.  We thus investigate their structure further,
    giving a formula for the multiplicity along every component of the
    jet schemes of a general reduced monomial hypersurface (that is,
    the case of a simple normal crossing divisor).
\end{abstract}

\thanks{I would like to thank Karen Smith for suggesting this problem,
and for many valuable discussions.  Also I am grateful to Hai Long Dao
for his inspiring idea.}

\maketitle

\section{Introduction}

Let $X$ be a scheme of finite type over an algebraically closed field
$k$ of characteristic zero.  An $m$-jet on $X$ is a $k$-morphism
$$\spec k[t]/(t^{m+1}) \rightarrow X.$$ We can think of such a jet as
the choice of a point $x$ in $X$, together with the choice of a
tangent direction, a second order tangent direction, and so forth up
to an $m^{th}$ order tangent direction at $x$.  The set of all
$m$-jets on $X$ forms a scheme in a natural way.  This is the $m^{th}$
jet scheme $\mathcal{J}_{m}(X)$.

The first serious study of arc spaces appears to have been undertaken
by Nash in the $60$'s \cite{Nash}.  Nash was interested in whether the
singularities of a variety $X$ could be reflected in its arc space.
Interest in jet schemes and arc spaces was recently rejuvenated by
Kontsevich's introduction of motivic integration in his $1995$ lecture
at Orsay \cite{Kontsevich}.  Since then the theory has been further
developed by Batyrev \cite{Batyrev1, Batyrev3}, Denef and Loeser
\cite{DenefLoeser1, DenefLoeser2, DenefLoeser4}, and others.  For good
surveys on motivic integration, see \cite{Blickle, Craw, DenefLoeser3,
Looijenga, Veys}.  Inspired by these developments but in a somewhat
different direction, Musta\c{t}\v{a} and his collaborators later
tied jet schemes to the study of singularities and invariants in 
birational geometry.  For results in this direction, see
\cite{FernexEinMustata, EinLazarsfeldMustata, EinMustata,
EinMustataYasuda, MustataLCI, MustataSingsPairs}.

However, there has been little study of the explicit scheme structure
of jet schemes, even in simple cases.  Goward and Smith
\cite{GowardSmith} began a study of the jet schemes of monomial
varieties by giving an explicit description of the irreducible
components of the jet schemes of monomial schemes, see Section
\ref{statement}.  In this paper, we continue this project by finding a 
formula for the multiplicity of these jet schemes along each of their 
irreducible components.

\section{Statement of main theorem}\label{statement}

First we recall how the set of all $m$-jets on an affine variety forms
a scheme.

Let $X \subseteq \mathbb{A}^{n}$ be a closed subscheme, say $X=\spec
k[x_{1},\ldots,x_{n}]/I$.  An $m$-jet of $X$ is equivalent to a
$k$-algebra homomorphism $$\phi:k[x_{1},\ldots,x_{n}]/I \rightarrow
k[t]/(t^{m+1}).$$
The map $\phi$ is clearly determined by the images of the generators
$x_{i}$ under $\phi$, say
$\phi(x_{i})=x_{i}^{(0)}+x_{i}^{(1)}t+\ldots+x_{i}^{(m)}t^{m}$.
However, there are constraints on the $x_{i}^{(j)}$'s coming from the
relations on the $x_{i}$'s.  Indeed, fixing a set of generators
$f_{1},\ldots,f_{d}$ for the ideal $I$, we must have $\phi(f_{k})=0$
for each $k=1,\ldots,d$, or that
\begin{equation}
    f_{k}(x_{1}^{(0)}+\ldots+x_{1}^{(m)}t^{m},\ldots,x_{n}^{(0)}
    +\ldots+x_{n}^{(m)}t^{m})=0 \hspace{0.2cm}\text{ in
    }k[t]/(t^{m+1}).  \tag{$\dagger$}
\end{equation} 
Rearranging ($\dagger$) so as to collect like terms of powers of $t$,
we have $$f_{k}^{(0)}+f_{k}^{(1)}t+\ldots+f_{k}^{(m)}t^{m}$$
where each $f_{k}^{(l)} \in
k[x_{1}^{(0)},\ldots,x_{1}^{(m)},\ldots,x_{n}^{(0)},\ldots,x_{n}^{(m)}]$.
Then $\mathcal{J}_{m}(X)$ is the subscheme of $\spec
k[x_{i}^{(j)}]=\mathcal{J}_{m}(\mathbb{A}^{n})$ defined by the
polynomials $f_{k}^{(l)}$, as $k$ ranges from $1$ to $d$ and $l$
ranges from $0$ to $m$.  When $X \subseteq \mathbb{A}^{n}$ is an
affine scheme, we let $J_{m}(X)$ denote the defining ideal of the jet
scheme $\mathcal{J}_{m}(X)$ as a subscheme of
$\mathcal{J}_{m}(\mathbb{A}^{n}) \cong \mathbb{A}^{n(m+1)}$.

For general properties of jet schemes, see \cite{Blickle, MustataLCI}.

\vspace{1cm}

Now let $X \subseteq \mathbb{A}^{n}$ be a reduced monomial hypersurface
defined by a single monomial $x_{1} \cdots x_{r}$.  Its jet scheme
$\mathcal{J}_{m}(X)$ is the closed subscheme of
$$\mathcal{J}_{m}(\mathbb{A}^{n})=\spec
k[x_{1}^{(0)},\ldots,x_{1}^{(m)},\ldots,x_{n}^{(0)},\ldots,x_{n}^{(m)}],$$
cut out by polynomials $g_{0},\ldots,g_{m}$ where $g_{k}$ is the
coefficient of $t^{k}$ in the product
$$\prod_{i=1}^{r}(x_{i}^{(0)}+x_{i}^{(1)}t+\ldots+x_{i}^{(m)}t^{m}).$$
In other words, the ideal $$J_{m}(X) \subseteq
R=k[x_{1}^{(0)},\ldots,x_{1}^{(m)},\ldots,x_{n}^{(0)},\ldots,x_{n}^{(m)}]$$
has generators of the form $$g_{k}\hspace{0.2cm}=\hspace{0.2cm} \sum
x_{1}^{(i_{1})}x_{2}^{(i_{2})}\cdots x_{r}^{(i_{r})}$$ where $\sum
i_{j}=k$ and $0 \leq i_{j} \leq m$.

\begin{theorem}[Goward-Smith]
    With $X$ as above, the minimal primes of $J_m(X)$ are of the form
    \begin{equation}
	P(m;t_1,\ldots,t_r)=(x_1^{(0)},\ldots,x_1^{(t_1-1)},\ldots,
	x_r^{(0)},\ldots,x_r^{(t_r-1)})  \tag{$\star$}
    \end{equation}
    where $0 \leq t_i \leq m+1$ and $\sum t_i=m+1$.  (Here, we adopt
    the convention that the value $t_i=0$ means the variable $x_i$
    does not appear at all.)
\end{theorem}    

In particular, $\mathcal{J}_{m}(X)$ is pure dimensional --- each
component has codimension $m+1$.  This gives a complete understanding
of the irreducible components of the reduced jet schemes of $X$, but
no more information about the scheme structure.  The following theorem
gives some insight into this scheme structure:

\begin{maintheorem}\label{Multiplicity}
    With notation as above, for $X=\spec
    k[x_{1},\ldots,x_{n}]/(x_{1}\cdots x_{r})$, the multiplicity of
    $\mathcal{J}_m(X)$ along $P(m;t_1, \ldots, t_r)$ is
    $$\frac{(m+1)!}{t_1!  \cdots t_r!}.$$
\end{maintheorem}

\section{Proof of main theorem}

We begin with two reductions.  Fix a minimal prime
$P=P(m;t_{1},\ldots,t_{r})$ of $J_{m}(X)$ as described in the theorem
above and let $R(m;t_1,\ldots,t_r)=R_{P}$.  We want to find $\ell
(R(m;t_1,\ldots,t_r)/J_m(X))$.

\begin{enumerate}
    \item Let $R'=k[x_i^{(j)}]$, $1 \leq i \leq r$, $j \geq 0$, be a
    polynomial ring in infinitely many indeterminates.  Also let
    $R'(m;t_1,\ldots,t_r)= R'_{P}$. Since length is unaffected by 
    completion, one can check then $\ell
    (R(m;t_1,\ldots,t_r)/J_m(X))=\ell (R'(m;t_1,\ldots,t_r)/J_m(X))$.
    So we may replace $R$ by $R'$.
    
    \item Suppose that one of the $t_{i}$, say $t_{r}$, is zero ---
    that is, suppose that $x_{r}$ does not appear in $P$ at all.
    Notice that we can rewrite the generators of $J_{m}(X)$ as
    $$g_{k}=\sum_{q=0}^{k} x_{r}^{(k-q)}h_{q} \hspace{1cm}
    \text{with}\hspace{1cm} h_{q}=\sum x_{1}^{(i_{1})}\cdots
    x_{r-1}^{(i_{r-1})}$$ where $\sum i_{j}=q$ and $0 \leq i_{j} \leq
    m$.  In other words,
    \begin{align*}
	g_{0} &= x_{r}^{(0)}h_{0}\\
	g_{1} &= x_{r}^{(1)}h_{0}+x_{r}^{(0)}h_{1}\\
	g_{2} &= x_{r}^{(2)}h_{0}+x_{r}^{(1)}h_{1}+x_{r}^{(0)}h_{2}\\
	&\text{etc.}
    \end{align*}
    Because $x_{r}^{(0)}$ is a unit in $R(m;t_{1},\ldots,t_{r-1},0)$,
    it follows that
    $$J_{m}(X)R(m;t_{1},\ldots,t_{r-1},0)=(h_{0},\ldots,h_{m})R(m;t_{1},\ldots,t_{r-1},0).$$
    This means that
    \begin{align*}
	\ell(R(m;t_{1},\ldots,t_{r-1},0)/J_{m}(X)) &=
	\ell(R(m;t_{1},\ldots,t_{r-1},0)/(h_{0},\ldots,h_{m}))\\
	&= \ell(R(m;t_{1},\ldots,t_{r-1})/J_{m}(X'))
    \end{align*}
    where $X'=\spec k[x_{1},\ldots,x_{n}]/(x_{1}\cdots x_{r-1})$.
    So we may assume all $t_{i}$ are positive.
\end{enumerate}

To prove the main theorem, we also need the following lemma:

\begin{lemma}\label{LengthLemma}
    Let $R$ be an arbitrary ring and $x_{1},\ldots,x_{r}$ nonunits in
    $R$.  If $x_1\cdots x_r$ is a nonzerodivisor on $R$, then $\ell
    \left(R/(x_1\cdots x_r) \right)=\sum_{i=1}^r \ell \left(R/(x_i)
    \right)$.
\end{lemma}
\begin{proof}
    The proof is a simple induction on $r$.
\end{proof}

\begin{proof}[Proof of Main Theorem]
    We proceed by induction on $m$.  A minimal prime of $J_0(X)$ has
    the form $(x_{k}^{(0)})$ and $\ell
    \left(k[x_{i}^{(j)}]_{(x_{k}^{(0)})}/(x_{1}^{(0)}\cdots
    x_{r}^{(0)})\right)= \ell
    \left(k[x_{i}^{(j)}]_{(x_{k}^{(0)})}/(x_{k}^{(0)})\right)=1$.
    This completes the $m=0$ case.  Now we assume that the assertion
    is true for $m-1$ and consider the multiplicity of
    $\mathcal{J}_m(X)$ along $P(m;t_1, \ldots, t_r)$.

    Note that
    \begin{align*}
	\height(J_m(X)R(m;t_1, \ldots, t_r)) &= \dim(R(m;t_1, \ldots,
	t_r))\\
	&= \sum_{i=1}^r t_i \hspace{1cm}\text{(by 
	$\star$)}\\
	&= m+1,
    \end{align*}
    which is the number of generators of $J_m(X)$.  So the elements
    $g_0, \ldots, g_m$ form a regular sequence on $R(m;t_1, \ldots,
    t_r)$.  Now let $S=R(m;t_1, \ldots, t_r)/(g_{1},\ldots,g_{m})$.
    Then $\ell (R(m;t_1, \ldots, t_r)/J_m(X))=\ell (S/(g_{0}))$.

    Since $g_0=x_1^{(0)} \cdots x_r^{(0)}$ is a nonzerodivisor on $S$,
    by the Lemma, $\ell \left(S/(x_1^{(0)} \cdots
    x_r^{(0)}) \right)$ $=\sum_{n=1}^r \ell \left(S/(x_n^{(0)}) \right)$.
    To know what $\ell \left(S/(x_n^{(0)}) \right)$ is, we need the
    following crucial result:

    \begin{claim}
    The ring $S/(x_{n}^{(0)})$ is isomorphic to
    $R(m-1;t_{1},\ldots,t_{n}-1,\ldots,t_{r})/J_{m-1}(X)$.
    \end{claim}
    \begin{proof}
	Modulo $x_{n}^{(0)}$, the generators $g_{k}$ $(1 \leq k \leq
	m)$ become $$\widetilde{g_{k}}=\sum
	x_{1}^{(i_{1})}x_{2}^{(i_{2})}\cdots x_{r}^{(i_{r})}$$ where
	$\sum i_{j}=k$, $0 \leq i_{j} \leq m$ if $j \neq n$ and $1
	\leq i_{n} \leq m$.

	By a change of variables, replacing $x_{n}^{(q)}$ by
	$x_{n}^{(q-1)}$ for all $q \geq 1$ and fixing the rest,
	$\widetilde{g_{k}}$ becomes $g_{k-1}$.  This tells us
	\begin{align*}
	    S/(x_{n}^{(0)}) &=
	    R(m;t_{1},\ldots,t_{r})/(\widetilde{g_{1}},\ldots,\widetilde{g_{m}},x_{n}^{(0)})\\
	    &\cong
	    R(m-1;t_{1},\ldots,t_{n}-1,\ldots,t_{r})/(g_{0},\ldots,g_{m-1})\\
	    &= R(m-1;t_{1},\ldots,t_{n}-1,\ldots,t_{r})/J_{m-1}(X).
	\end{align*}
    \renewcommand{\qed}{}
    \end{proof}

    Therefore, by the induction hypothesis, we have
    \begin{align*}
	\ell(R(m;t_1, \ldots, t_r)/J_m(X)) &= \sum_{n=1}^{r}
	\frac{m!}{t_{1}!\cdots (t_{n}-1)!\cdots t_{r}!}\\
	&= \frac{m!(t_1+\ldots+t_r)}{t_1!\cdots t_r!}\\
	&= \frac{(m+1)!}{t_1!\cdots t_r!}.
    \end{align*}
\end{proof}

\bibliographystyle{plain}
\bibliography{biblio}   % Use the BibTeX file ``biblio.bib''.

\end{document}